\def\marker{\>\hbox{${\vcenter{\vbox{
    \hrule height 0.4pt\hbox{\vrule width 0.4pt height 6pt
    \kern6pt\vrule width 0.4pt}\hrule height 0.4pt}}}$}\>}

\documentclass[12pt]{article}
\usepackage{fullpage}
\usepackage{amsmath,amsfonts,amssymb,amsthm,graphics,amscd,MnSymbol}
\usepackage{graphicx}

\usepackage[
pdfauthor={},
pdftitle={Generalized Petersen graphs are $(1,3)$-choosable},
pdfstartview=XYZ,
bookmarks=true,
colorlinks=true,
linkcolor=blue,
urlcolor=blue,
citecolor=blue,
bookmarks=false,
linktocpage=true,
hyperindex=true
]{hyperref}

\newtheorem {Theorem}  {Theorem}[section]
\newtheorem {Lemma}[Theorem]{Lemma}

\newtheorem {Corollary}[Theorem]{Corollary}

\theoremstyle{definition}
\newtheorem{Definition}[Theorem]{Definition}

\newtheorem {Observation}[Theorem]{Observation}
\newenvironment {Proof} {\noindent {\bf Proof.}}{\quad $\square$\par\vspace{3mm}}

\def\mad{{\rm mad}}

\date{\today}
\title{\bf Generalized Petersen graphs are $(1,3)$-choosable}

\author{Yunfang Tang\thanks{
Department of Mathematics, China Jiliang University, Hangzhou 310018, P.R. China. Grant number:
  NSF11701543. Email: tangyunfang8530@cjlu.edu.cn.}
\and Yuting Yao \thanks{
Department of Mathematics, China Jiliang University, Hangzhou 310018, P.R. China. Email: s22080701014@cjlu.edu.cn.}
\\[0.2cm]
      }

\begin{document}

\maketitle

\vspace{-2pc}

\begin{abstract}
A total weighting of a graph $G$ is a mapping $\phi$
that assigns a weight to each vertex and each edge of $G$.
The vertex-sum of $v \in V(G)$ with respect to $\phi$
is $S_{\phi}(v)=\sum_{e\in E(v)}\phi(e)+\phi(v)$.
A total weighting is proper if adjacent vertices
have distinct vertex-sums.
A graph $G=(V,E)$ is called $(k,k')$-choosable if the
following is true: If each vertex $x$ is assigned a set $L(x)$ of $k$ real numbers,
and  each edge $e$ is assigned a set $L(e)$ of $k'$ real numbers,
then there is a proper total weighting $\phi$ with $\phi(y)\in L(y)$
for any $y \in V \cup E$.
In this paper, we prove that the generalized Petersen graphs are $(1,3)$-choosable.

\end{abstract}

\section{Introduction}

A \textit{total weighting} of a graph $G$ is a mapping $\phi: V(G)\cup E(G)\rightarrow \mathbb{R}$.
The {\em vertex-sum} of $v$ with respect to $\phi$ is
$S_{\phi}(v)=\sum_{e\in E(v)}\phi(e)+\phi(v)$, where $E(v)$ is the set of edges incident to $v$.
A total weighting $\phi$
with $\phi(v)=0$ for all vertices $v$ is also called an {\em edge weighting}.
A total weighting $\phi$ is \textit{ proper } if for any edge $uv$ of $G$,
$S_{\phi}(u) \ne S_{\phi}(v)$.

Proper edge weighting of graphs was first studied by Karo\'nski, {\L}uczak and Thomason \cite{KLT2004}.  They conjectured that every graph with no isolated edge has a proper edge-weighting using weights 1,2, and 3.
This conjecture is called the 1-2-3 Conjecture, and attracted considerable attention, and it finally proved by Ralph Keusch in  \cite{K2023arxiv}.

The list version of edge weighting of graphs was introduced by Bartnicki, Grytczuk and
Niwczyk in \cite{BGN2009}. The list version of total weighting of graphs was introduced independently
by Przyby{\l}o and Wo\'{z}niak \cite{PW2011} and by Wong and Zhu in \cite{WZ2011}.
Suppose $\psi:V(G)\cup E(G) \longrightarrow\{1,2,\ldots\}$
is a mapping which assigns to each vertex and each edge of $G$ a positive integer.
A $\psi$-list assignment of $G$ is a mapping $L$ which assigns to $z\in V(G)\cup E(G)$
a set $L(z)$ of $\psi(z)$ real numbers. Given a total list assignment $L$,
a \textit{proper $L$-total weighting} is a proper total weighting $\phi$ with
$\phi(z)\in L(z)$ for all $z\in V(G)\cup E(G)$.
We say $G$ is \textit{total weight $\psi$-choosable} if for any
$\psi$-list assignment $L$, there is a proper $L$-total weighting of $G$.
We say $G$ is $(k,k')$-choosable if $G$ is $\psi$-total weight choosable,
where $\psi(v)=k$ for $v\in V(G)$
and $\psi(e)=k'$ for $e\in E(G)$.

It was conjectured in \cite{WZ2011} that every graph with
no isolated edges is $(1,3)$-choosable and every graph is $(2,2)$-choosable, which are the list version of 1-2-3 Conjecture and list version of 1-2 Conjecture, respectively. Wong and Zhu \cite{WZ2016} proved  that
every graph is $(2,3)$-choosable. Cao \cite{C2021} proved that every graph with no isolated edges is $(1,17)$-choosable, and  Zhu \cite{Z2022} proved that every graph with no isolated edges is $(1,5)$-choosable.
But the $(1,3)$-choosable Conjecture and the $(2,2)$-choosable Conjecture remain open. Indeed, it
 is unknown whether there is a constant $k$
such that every graph is $(k,2)$-choosable.

Some special graphs are shown to be $(1,3)$-choosable, such as complete graphs,
complete bipartite graphs, trees (these graphs are shown in \cite{BGN2009}), connected $2$-degenerate non-bipartite graphs other than $K_{2}$ \cite{WZ2015},
graphs with maximum average degree less than $\frac{11}{4}$ \cite{LWZ2018}, and Halin graphs \cite{LWZ2021}.

The  \textit{generalized Petersen graph}  $ P(n,t)$ is the graph   with vertex set $V(G)=\{u_{i},v_{i}: 1\leq i \leq n\}$ and edge set $E(G)=\{u_{i}u_{i+1},u_{i}v_{i},v_{i}v_{i+t}\}$, where the indices are modulo $n$.

This paper proves that every generalized Petersen graph  is $(1,3)$-choosable.

\section{Preliminaries}
Assume $G$ is a simple graph, where $V(G)=\{v_1,v_2,\ldots,v_n\}$ and $E(G)=\{e_1,e_2,\ldots,e_m\}$. Let $T(G)=V(G)\cup E(G)$. We assign a variable $x_v$ to each vertex, and $x_e$ to each edge.  For each vertex $v$ of $G$, let $S_u=\sum_{e\in E(u)}x_e+x_u$. Fix an arbitrary orientation $D$ of $G$. Let
\begin{eqnarray*}
P(G)=P(x_{v_1},x_{v_2},\ldots, x_{v_n},x_{e_1},x_{e_2},\ldots,x_{e_m})=\prod_{u\in V(G)}P(G,u),
\end{eqnarray*}
where
\begin{eqnarray*}
P(G,u)=\prod_{v\in N^{+}_D(u)}(S_u-S_v)
\end{eqnarray*}

It follows from the definition that
a mapping $\phi:T(G)\rightarrow \emph{R}$ is a proper total weighting of $G$ if and only if  $P(\phi) \ne 0$, where $P(\phi)$ is the evaluation of $P$ at $x_z=\phi(z) $
for   $z\in T(G)$.

Assume $X$ is a subset of $V(G)$ and $H=G[X]$ is the subgraph of $G$ induced by $X$.    Denote by $T(G)\setminus T(H)=\left(V(G)\setminus V(H) \right)\cup\left(E(G)\setminus E(H)\right)$.
\begin{Definition}
An induced subgraph $H$ of $G$ is called a \textit{$(k,k')$-choosable induced subgraph } in $G$ if for any $(k,k')$-total list assignment $L$ of $G$ and any $L$-total weighting $\psi:T(G)\setminus T(H)\rightarrow \emph{R}$, there
is an $L$-total weighting $\phi$ of $G$ such that
\begin{itemize}
\item[(i)] $\phi\mid_{T(G)\setminus T(H)}=\psi$;
\item[(ii)] $\sum_{z\in E_G(u)\cup\{u\}}\phi(z)\neq \sum_{z\in E_G(v)\cup\{v\}}\phi(z)$ for every edge $uv\in E(G)$ with $u\in X$.
\end{itemize}
\end{Definition}
Clearly, if $G$ has a $(k,k')$-choosable induced subgraph $H$ such that $T(G)\setminus T(H)$ has an $L$-total weighting $\psi$ satisfying $\sum_{z\in E_G(u)\cup{\{u\}}}\psi(z)\neq \sum_{z\in E_G(v)\cup{\{v\}}}\psi(z)$ for any two adjacent vertices $u,v\in V(G)\setminus V(H)$, then $G$ is $(k,k')$-choosable. Especially, we have straightforward result in the following:

\begin{Lemma} \label{key1}
Suppose that $G$ has a $(k,k')$-choosable induced subgraph $H$. Let $G'=G[E(G)\setminus E(H)]$.  If $G'$ is $(k,k')$-choosable, then $G$ is $(k,k')$-choosable.
\end{Lemma}
\begin{Proof}
If $G'$ is $(k,k')$-choosable, then $G'$ has an $L$-total weighting $\psi$ satisfying $$\sum_{z\in E_{G'}(u)\cup{
\{u\}}}\psi(z)\neq \sum_{z\in E_{G'}(v)\cup{\{v\}}}\psi(z)$$ for any two adjacent vertices $u,v\in V(G)\setminus V(H)\subseteq V(G')$. Note that $E_G(u)=E_{G'}(u)$ for each vertex $u\in V(G)\setminus V(H)$. Then $T(G)\setminus T(H)$ also has an $L$-total weighting $\psi$ satisfying $\sum_{z\in E_G(u)\cup{\{u\}}}\psi(z)\neq \sum_{z\in E_G(v)\cup{\{v\}}}\psi(z)$ for any two adjacent vertices $u,v\in V(G)\setminus V(H)$. Thus the result follows.
\end{Proof}

Let $L$ be a $(k,k')$-total list assignment of $G$ and $\psi$ be an $L$-total weighting on $T(G)\setminus T(H)$.
Assign a variable $x_z$ to each vertex or edge $z\in T(G)$ satisfying $x_z=\psi(z)$ is a constant for $z\in T(G)\setminus T(H)$. Fix an orientation $D$ on $E(G)$ such that every edge of $\cup_{u\in V(H)}E_G(u)\setminus E(H)$ is oriented away from $H$. Then
we define a polynomial $Q_G(H)$ by
\begin{eqnarray*}
Q_G(H)=Q(\{x_z:z\in T(H)\})=\prod_{u\in V(H)}Q_G(H,u)=\prod_{u\in V(H)}\prod_{v\in N_{D}^{+}(u)}\left(\sum_{z\in E_G(u)\cup\{u\}}x_z-\sum_{z\in E_G(v)\cup\{v\}}x_z\right)
\end{eqnarray*}
Note that $H$ is a $(k,k')$-choosable induced subgraph in $G$ if and only if $L$ and $\psi$ are defined as above, there is an $\phi(z)\in L(z)$ for each $z\in T(H)$ such that $Q(\{\phi(z):z\in T(H)\})\ne 0$.

We will apply the following Combinatorial Nullstellsatz to establish a sufficient condition for getting a $(k,k')$-choosable induced subgraph.

\begin{Lemma} (Combinatorial Nullstellensatz {\rm\cite{A1999}}) \label{CN} Let $\mathbb{F}$ be an arbitrary field, and let $P = P(x_{1},x_{2},\ldots,x_{n})$ be a polynomial in $\mathbb{F}[x_{1},x_{2},\ldots,x_{n}]$. Suppose the degree $deg(P)$ of $P$ equals $\sum_{i=1}^{n}k_{i}$, where each $k_{i}$ is a non-negative integer, and suppose the coefficient of $x_{1}^{k_{1}}\cdots$$x_{n}^{k_{n}}$ in $P$ is non-zero. Then if $S_{1},\ldots,S_{n}$ are subsets of $\mathbb{F}$ with $\mid S_{i}\mid>k_{i}$, there are $s_{1}\in S_{1},\ldots,s_{n}\in S_{n}$ so that $P(s_{1},s_{2},\ldots,s_{n})\neq 0$.
\end{Lemma}

The goal of this paper is to prove that $H$ is a $(1,3)$-choosable induced subgraph of $G$ and moreover, $x_z=\psi(z)$ is a constant $z\in T(G)\setminus T(H)$. Thus we restrict to
monomials of the form $\prod_{e\in E(H)}x_e^{i_e}$ by omitting the variables $x_v$ for $v\in V(H)$ and vanish all monomials $\prod_{z\in T(H)}x_z^{i_z}$ with $\sum_{z\in T(H)}i_z<\deg(Q_G(H))$. For this purpose, we  will consider the following polynomial:
\begin{eqnarray*}
P_G(H)&=&P(\{x_e:e\in E(H)\})=\prod_{u\in V(H)}P_G(H,u)\\
&=&\prod_{u\in V(H)}\left(\left(\sum_{e\in E_H(u)}x_e\right)^{|E_G(u)\setminus E_H(u)|}
\prod_{v\in N_{D_H}^{+}(u)}\left(\sum_{e\in E_G(u)}x_e-\sum_{e\in E_G(v)}x_e\right)\right),
\end{eqnarray*}
where $D_H$ is the restriction of $D$ on $E(H)$.

By Lemma~\ref{CN}, the following result is straightforward observation.

\begin{Corollary}\label{key2}
Let $H$ be an induced subgraph of a graph $G$. If $P_G(H)$ has a monomial $\prod_{e\in E(H)}x_e^{i_e}$ with non-zero coefficient, then $H$ is a $(1,3)$-choosable induced subgraph in $G$, where $|i_e|\le 2$ for each edge $e\in E(H)$.
\end{Corollary}

\begin{Definition}\label{NJ} Let $X=\{x_{1},x_{2},\ldots,x_{n}\}$ and $Y=\{x_{i_{1}},\ldots,x_{i_{m}}\} \subseteq X$. We can think of every polynomial in $\mathbb{F}[X]$ as a polynomial in $\mathbb{T}[Y]$ with $\mathbb{T}=\mathbb{F}[X \backslash Y] $. For $J = x_{i_{1}}^{j^{1}}\cdots x_{i_{m}}^{j_{m}}$, we define a mapping $\eta_{J}:\mathbb{F}[X]\rightarrow\mathbb{T}$ by $\eta_{J}[P]$ = the coefficient of J in P $\in \mathbb{T}[Y]$, for any $P \in \mathbb{F}[X]$.
\end{Definition}

\begin{Lemma} \label{use}
Let $P(x_{1},x_{2},\ldots,x_{n})$ be a polynomial of degree $m$ with $n$ variables, where $m,n\geq1$.  Let
\begin{displaymath}
P(x_{1},x_{2},\ldots,x_{n})=\prod_{k=1}^{s}P_{k}, \qquad J=\prod_{k=1}^{s}J_{k},
\end{displaymath}
where $s$ is a positive integer and $1\leq s\leq m $. For $k \geq 2$,
If $P_{k}$ does not contain any variable in $J_{1},J_{2},\ldots,J_{k-1}$,  $deg(\prod_{k=1}^{t}J_{k})=deg(\prod_{k=1}^{t}P_{k})$, and $deg(J_{i})=deg(P_{i})$ for $t+1\leq i\leq s$, where $t$ is a positive integer and $1\leq t\leq s$, then
$$\eta_{J}[P] =\eta_{\prod_{k=1}^{t}J_{k}}[\prod_{k=1}^{t}P_{k}]\cdot\prod_{k=t+1}^{s}\eta_{J_{k}}[P_{k}].$$
\end{Lemma}

\begin{Proof}
\begin{eqnarray*}
\eta_{J}[P] &=&\eta_{J_{t+1}\cdots J_{s}}[\eta_{\prod_{k=1}^{t}J_{k}}[\prod_{k=1}^{t}P_{k}]\cdot \prod_{k=t+1}^{s}P_{k}]\\
&=&\eta_{\prod_{k=1}^{t}J_{k}}[\prod_{k=1}^{t}P_{k}]\cdot\eta_{J_{t+1}\cdots J_{s}}[\prod_{k=t+1}^{n}P_{k}]\\
&=&\eta_{\prod_{k=1}^{t}J_{k}}[\prod_{k=1}^{t}P_{k}]\cdot\prod_{k=t+1}^{s}\eta_{J_{k}}[P_{k}].
\end{eqnarray*}
\end{Proof}

With the help of MATHEMATICA, the following results we will use in Lemma~\ref{L1} are directly obtained.
\begin{Observation}\label{Ob2}
For any positive integer $s$ and $t$, and $s< t$, denote by $J^*_{s,t}=\prod_{i=s}^{t}y_{i}^{2}$,
$F^*_{s,t}=\prod_{i=s}^{t}(y_{i}+y_{i+1})(y_{i}-y_{i+2})$, then $\eta_{J^*_{s,t}}[F^*_{s,t}]=1$.
\end{Observation}

\begin{Theorem}{\rm\cite{BGN2009}}\label{tree}
If  $G$ is a forest without isolated edges, then $G$ is $(1,3)$-choosable.
\end{Theorem}

\begin{Theorem}{\rm\cite{LWZ2018}}\label{mad}
Graphs with maximum average degree less than $\frac{11}{4}$ are $(1,3)$-choosable.
\end{Theorem}

\begin{Theorem}{\rm\cite{WZ2015}}\label{2d}
Every connected $2$-degenerate non-bipartite graph other than $K_{2}$ is $(1,3)$-choosable.
\end{Theorem}

\section{Proof of Theorem}
\begin{Lemma}\label{L1}
Let $G$ be a cubic graph, and $H$ be an induced subgraph of $G$. If $H$ is one of the following configurations:
\begin{itemize}
\item [(i)] $H$ is a graph obtained by two adjacent cycles  with a common edge, where one cycle is an $n$-cycle $u_{1}u_{2}\cdots u_{n}u_1$ and the other one is a $4$-cycle $u_1u_2v_2v_1u_1$;
\item [(ii)] $H$ is a graph obtained by two adjacent $5$-cycles with a common edge, which is also called diamond $D_8$;
\item [(iii)] $H$ is a graph obtained by by two adjacent cycles  with two common edges, where one cycle is an $n$-cycle $u_{1}u_{2}\cdots u_{n}u_1$ and the other one is a $8$-cycle $u_{1}u_{2}v_{2}v_{t+2}u_{t+2}u_{t+1}v_{t+1}v_{1}u_1$.
\end{itemize}
then $H$ is a $(1,3)$-choosable induced subgraph in $G$.
\end{Lemma}

\begin{figure}[h]
\centering
\includegraphics[scale=0.8]{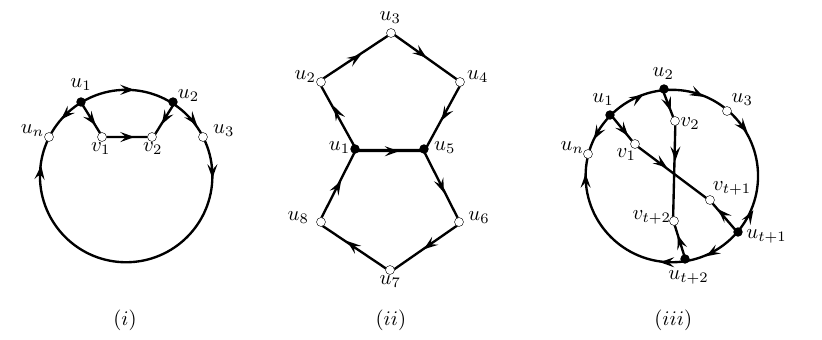}\\
\small{Fig. 1. (i)-(iii)}
\end{figure}

\begin{Proof}
In the proof of this Lemma, the coefficients of all the polynomials we use are obtained by MATHEMATICA.

(i) Assign variables $y_{i}$ to $u_{i}u_{i+1}$ for $i=1,2,\ldots,n$; $z_{i}$ to $u_{i}v_{i}$ for $i=1,2$, and $w$ to $v_{1}v_{2}$, respectively. Fix an orientation of $H$ such that every edge $u_{i}u_{i+1}$ is oriented toward $u_{i+1}$ for $i\in[1,n-1]$, every edge $u_{i}v_{i}$ is oriented toward $v_{i}$ for $i\in\{1,2\}$, edge $u_{1}u_{n}$ is oriented toward $u_{n}$ and edge $v_{1}v_{2}$ is oriented toward $v_{2}$ (see Figure $1(i)$).
Let $V(H_0)=\{u_{i},v_{j}:i,j\in\{1,2\}\}$. Then
\begin{eqnarray*}
P_G(H)&=&\prod_{v\in V(H_0)}P_G(H,v)\cdot\prod_{i=3}^{i=n}P_G(H,u_i)\\
&=&[(y_{n}+y_{1}-w)(y_{n}+z_{1}-y_{2}-z_{2})(y_{1}+y_{2}-w)(y_{1}+z_{2}-y_{3})\\
&&(y_{1}+z_{1}-y_{n-1})(z_{1}+w)(z_{1}-z_{2})(z_{2}+w)]\cdot F^*_{2,n-2}\cdot(y_{n-1}+y_{n}).
\end{eqnarray*}

Choose $$J=z_{1}^{2}z_{2}w^{2}\prod_{i=1}^{n-1}y_{i}^{2}=(z_{1}wy_1)^{2}z_{2}\cdot J^*_{2,n-2}\cdot y_{n-1}^{2}.$$
Then by Lemma~\ref{use} and Observation~\ref{Ob2}, we obtain that
\begin{eqnarray*}
\eta_{J}[P_G(H)]&=&\eta_{y_{n-1}^{2}}\{\eta_{(z_{1}wy_1)^{2}z_{2}}[(y_{n}+y_{1}-w)(y_{n}+z_{1}-y_{2}-z_{2})(y_{1}+y_{2}-w)(y_{1}+z_{2}-y_{3})\\
&&(y_{1}+z_{1}-y_{n-1})(z_{1}+w)(z_{1}-z_{2})(z_{2}+w)]\cdot\eta_{J^*_{2,n-2}}[F^*_{2,n-2}]\cdot(y_{n-1}+y_{n})\}\\
&=&\eta_{y_{n-1}^{2}}[-3y_{n-1}\cdot (y_{n-1}+y_{n})]\\
&=&-3 \neq0.
\end{eqnarray*}
Thus the result $(i)$ follows by Corollary~\ref{key2}.

(ii) We may assume that $H=C_8+u_1u_5$, where $C_8=u_1u_2\cdots u_8u_1$. For convenience, set $u_9=u_0$.
Assign variables $y_{i}$ to $u_{i}u_{i+1}$ for $i\in[1,8]$, and $z_{1}$ to $u_{1}u_{5}$. Fix an orientation of $H$ such that every edge $u_{i}u_{i+1}$ is oriented toward $u_{i+1}$ for $i\in[1,8]$, and edge $u_{1}u_{5}$ is oriented toward $u_{5}$ (see Figure $1(ii)$). Then
\begin{eqnarray*}
P_G(H)&=&\prod_{v\in V(C_8)}P_G(H,v)=(y_{8}+y_{1}-y_{4}-y_{5})(y_{8}+z_{1}-y_{2})(y_{3}-y_{5}-z_{1})(y_{4}+z_{1}-y_{6})\\
&&(y_{7}-y_{1}-z_{1})\cdot\prod_{i=1,2,3,5,6,7}(y_{i}+y_{i+1})\cdot\prod_{i=1,2,5,6}(y_{i}-y_{i+2}).
\end{eqnarray*}
Let $J=z_{1}^{2}\prod_{i=1}^{7}y_{i}^{2}/ y_{4}$, we obtain that $\eta_J[P_G(H)]=-8\neq0$. Thus the result follows by Corollary~\ref{key2}.

(iii) Assign variables $y_{i}$, $z_{j}$ and $w_{k}$ to $u_{i}u_{i+1}$, $u_{j}v_{j}$, and $v_{k}v_{k+t}$, respectively, for $i=\{1,2,\ldots,n\}; j=\{1,2,t+1,t+2\}; k=\{1,2\}$, where the subscripts of $u_{i}$ are taken modulo $n$.
Fix an orientation of $H$ such that every edge $u_{i}u_{i+1}$ is oriented toward $u_{i+1}$ for $i\in[1,t-1]\cup[t+1,n-1]$, edge $u_{j}v_{j}$ is oriented toward $v_{j}$ for $j\in\{1,2,t+1,t+2\}$, edge $u_{1}u_{n}$ is oriented toward $u_{n}$, edge $u_{t+1}u_{t}$ is oriented toward $u_{t}$, and edge $v_{k}v_{k+t}$ is oriented toward $v_{k+t}$ for $k\in\{1,2\}$ (see Figure $1(iii)$). Note that $V(C_8)=\{u_{1},u_{2},v_{2},v_{t+2},u_{t+2},u_{t+1},v_{t+1},v_{1}\}$. Then
\begin{eqnarray*}
P_G(H)=\prod_{i=1}^{t-1}\prod_{i=t+2}^{n-1}P_{i},
\end{eqnarray*}
where
\begin{eqnarray*}
P_{1}&=&\prod_{v\in V(C_8)}P_G(H,v)=
(y_n+z_1-y_2-z_2)(y_n+y_1-w_1)(y_1+z_1-y_{n-1})(y_1+z_2-y_3)\\
&&(y_1+y_2-w_2)(y_t+y_{t+1}-w_1)(y_{t}+z_{t+1}-y_{t+2}-z_{t+2})(y_{t+1}+y_{t+2}-w_2)(y_{t+1}+z_{t+1}-y_{t-1})\\
&&(y_{t+1}+z_{t+2}-y_{t+3})(z_{1}-z_{t+1})(z_1+w_1)(z_{2}-z_{t+2})(z_2+w_2)(z_{t+1}+w_1)(z_{t+2}+w_2),\\
P_{i}&=&P_G(H,u_{i+1})=\begin{cases}
(y_{i}+y_{i+1})(y_{i}-y_{i+2}) ~& \mbox{if $i\in[2,t-2]\bigcup[t+2,n-2]$},\\
(y_{i}+y_{i+1}) ~& \mbox{if $i\in \{t-1, n-1\}$}.
\end{cases}
\end{eqnarray*}

Choose $J=\prod_{i=1}^{t-1}\prod_{i=t+2}^{n-1}J_{i}$, where
\begin{eqnarray*}J_{i}&=&
\begin{cases}
y_{1}^{2}z_{1}^{2}z_{2}^{2}z_{t+1}^{2}z_{t+2}^{2}w_{1}^{2}w_{2}^2~& \mbox{if $i=1$},\\
y_{i}^{2}~& \mbox{if $i\in[2,t-2]\bigcup[t+2,n-2]$},\\
y_{t-1}^{2}y_{t} ~&\mbox{if $i=t-1$, and $t\neq 2~( mod ~3)$},\\
y_{t-1}y^{2}_{t} ~&\mbox{if $i=t-1$, and $t= 2~( mod ~6)$},\\
y_{i}y_{i+1} ~&\mbox{if $i\in\{t-1,n-1\}$, and $t= 5~( mod~ 6)$},\\
y_{n-1}~& \mbox{if $i=n-1$, and $t\neq 5~(mod~ 6)$}.
\end{cases}
\end{eqnarray*}

In the following, we will determine $\eta_J[P_G(H)]$ by some procedures.
For convenience, denoted by
\begin{eqnarray*}
Q^*=-y_n^2-y_{n-1}y_{t-1}-y_ny_{t-1}-y_{n-1}y_t+y_ny_t-y_t^2-2y_ny_{t+2}+2y_ty_{t+2}-y_{t+2}^2,
\end{eqnarray*}
$P'_{1}=P_{1}$, and $P'_{i}=\eta_{J_{i-1}}[P'_{i-1}]\cdot P_{i}$
for $i\in[2,t-1]$.

With the help of MATHEMATICA, we determine that $\eta_{J_{i}}[P'_{i}]$ for $i\in[1,t-1]$ in the following.
For $i=3k+1$,
\begin{eqnarray*}
&&\eta_{J_{1}}[P'_{1}]=\eta_{J_{1}}[P_{1}]\\
&=&-y_{2}^{2}+y_2(2y_n-2y_t+y_{t+2}-y_{t+3})-y_3(y_{t+2}+y_{t+3})+Q^*,\\
&&\eta_{J_{3k+1}}[P'_{3k+1}]=\eta_{J_{3k+1}}[\eta_{J_{3k}}[P'_{3k}]\cdot P_{3k+1}]\\
&=&-y_{3k+2}^{2}+(-1)^{3k+1}[y_{3k+2}(-2y_n+2y_t-y_{t+2}+y_{t+3})+y_{3k+3}(y_{t+2}+y_{t+3})]+Q^*;
\end{eqnarray*}
for $i=3k+2$,
\begin{eqnarray*}
&&\eta_{J_{2}}[P'_{2}]=\eta_{J_{2}}[\eta_{J_{1}}[P'_{1}]\cdot P_2]\\
&=&y_{3}y_4+2y_{3}(y_n-y_t-y_{t+3})+y_4[2(y_{t}-y_{n})-y_{t+2}+y_{t+3}]+Q^*,\\
&&\eta_{J_{3k+2}}[P'_{3k+2}]=\eta_{J_{3k+2}}[\eta_{J_{3k+1}}[P'_{3k+1}]\cdot P_{3k+2}]\\
&=&y_{3k+3}y_{3k+4}+(-1)^{3k+2}[2 y_{3k+3}(y_n-y_t-y_{t+3})+y_{3k+4}(2y_t-2y_n-y_{t+2}+y_{t+3})]+Q^*;
\end{eqnarray*}
and for $i=3k+3$,
\begin{eqnarray*}
&&\eta_{J_{3}}[P'_{3}]=\eta_{J_{3}}[\eta_{J_{2}}[P'_{2}]\cdot P_3]\\
&=&y_4^2-y_4(y_{t+2}+y_{t+3}+y_5)-2y_5(y_n-y_t-y_{t+3})+Q^*,\\
&&\eta_{J_{3k+3}}[P'_{3k+3}]=\eta_{J_{3k+3}}[\eta_{J_{3k+2}}[P'_{3k+2}]\cdot P_{3k+3}]\\
&=&y_{3k+4}^{2}-y_{3k+4}y_{3k+5}+(-1)^{3k+3}[y_{3k+4}(y_{t+2}+y_{t+3})+2y_{3k+5}(y_n-y_t-y_{t+3})]+Q^*.
\end{eqnarray*}

Set $i=t-2$. Then by the above computations, \begin{eqnarray*}
&&\eta_{J_{t-2}}[P'_{t-2}]\\
&=&\begin{cases}
-y_{t-1}^{2}+(-1)^{t}[y_{t-1}(-2y_n+2y_t-y_{t+2}+y_{t+3})+y_{t}(y_{t+2}+y_{t+3})]+Q^* &~\mbox{if $t=0~(mod~ 3)$},\\
y_{t-1}y_{t}+(-1)^{t}[2y_{t-1}(y_n-y_t-y_{t+3})+y_{t}(2y_t-2y_n-y_{t+2}+y_{t+3})]+Q^* &~\mbox{if $t=1~(mod~ 3)$},\\
y_{t-1}^{2}-y_{t-1}y_{t}+(-1)^{t-2}[y_{t-1}(y_{t+2}+y_{t+3})+2y_{t}(y_n-y_t-y_{t+3})]+Q^* &~\mbox{if $t=2~(mod~ 3)$}.
\end{cases}
\end{eqnarray*}
Recall that $P'_{t-1}=\eta_{J_{t-2}}[P'_{t-2}]\cdot(y_{t-1}+y_t)$.
Then \begin{eqnarray*}
\eta_{J_{t-1}}[P'_{t-1}]=\begin{cases}
-1+2(-1)^{t} &~\mbox{if $t=0~(mod~ 3)$},\\
1-2(-1)^{t} &~\mbox{if $t=1~(mod~ 3)$},\\
-4 &~\mbox{if $t=2~(mod~ 6)$},\\
-2y_{n-1}-2y_{n}+y_{t+2}+y_{t+3}&~\mbox{if $t=5~(mod~ 6)$}.
\end{cases}
\end{eqnarray*}

Suppose that $t\ne 5~(mod ~6)$. It's easily seen that $\deg(\prod_{i=1}^{t-1}J_{i})=\deg(\prod_{i=1}^{t-1}P_{i})$, $\deg(J_{i})=\deg(P_{i})$ and $P_{i}$ does not contain any variable in $J_{1},J_{2},\ldots,J_{i-1}$ for $i\in[t+2, n-1]$. Then by Lemma~\ref{use} and Observation~\ref{Ob2},
we have
\begin{eqnarray*}
\eta_{J}[P_{G}(H)]
&=&\eta_{\prod_{k=1}^{t-1}\prod_{k=t+2}^{n-1}J_{k}}[\prod_{k=1}^{t-1}\prod_{k=t+2}^{n-1}P_{k}]\\
&=&\eta_{\prod_{k=1}^{t-1}J_{k}}[\prod_{k=1}^{t-1}P_{k}]\cdot \prod_{k=t+2}^{n-1}\eta_{J_{k}}[P_{k}]\\
&=&\eta_{J_{t-1}}[P'_{t-1}]\cdot\eta_{J^*_{t+2,n-2}}[F^*_{t+2,n-2}]\cdot\eta_{J_{n-1}}[P_{n-1}]\\
&=&
\begin{cases}
-1+2(-1)^{t} &~\mbox{if $t=0~(mod~ 3)$},\\
1-2(-1)^{t} &~\mbox{if $t=1~(mod~ 3)$},\\
-4 &~\mbox{if $t=2~(mod~ 6)$}.
\end{cases}
\end{eqnarray*}

In the following, we will consider $t= 5~(mod ~6)$.
For convenience, let $Q_0^*=-2y_{n-1}-2y_{n}$, $P'_{t+1}=\eta_{J_{t-1}}[P'_{t-1}]$, and $P'_{i}=\eta_{J_{i-1}}[P'_{i-1}]\cdot P_{i}$, where $i\in[t+2,n-1]$.
With the help of MATHEMATICA, for $s\in[2,n-t+2]$,
\begin{eqnarray*}
\eta_{J_{t+s}}[P'_{t+s}]=
\begin{cases}
Q_0^*+(-1)^{s+1}(y_{t+s+1}+y_{t+s+2})~& \mbox{if $s=1~(mod~ 3)$},\\
Q_0^*+(-1)^{s}(2y_{t+s+1}-y_{t+s+2})~& \mbox{if $s=2~(mod~ 3)$},\\
Q_0^*+(-1)^{s+1}(y_{t+s+1}-2y_{t+s+2})~& \mbox{if $s=0~(mod~ 3)$}.
\end{cases}
\end{eqnarray*}
Let $s=n-t-2$. Then $n-2=t+s=3k+2+s$, we have
\begin{eqnarray*}
\eta_{J_{n-2}}[P'_{n-2}]=
\begin{cases}
[(-1)^{n}-2](y_{n-1}+y_n)~& \mbox{if $n=2~(mod~ 3)$},\\
[(-1)^{n-1}-1]2y_{n-1}-[2+(-1)^{n-1}]y_{n}~& \mbox{if $n=0~(mod~ 3)$},\\
[(-1)^{n}-2]y_{n-1}-[1+(-1)^{n}]2y_{n}~& \mbox{if $n=1~(mod~ 3)$}.
\end{cases}
\end{eqnarray*}
Note that $P'_{n-1}=\eta_{J_{n-2}}[P'_{n-2}]\cdot(y_{n-1}+y_n)$, and $J_{n-1}=y_{n-1}y_{n}$, then
\begin{eqnarray*}
\eta_{J}[P_G(H)]=\eta_{J_{n-1}}[P'_{n-1}]&=&
\begin{cases}
(-1)^{n-1}-4~& \mbox{if $n\neq2~(mod ~3)$}\\
2(-1)^{n}-4~& \mbox{if $n=2~(mod~ 3)$}\\
\end{cases}
\end{eqnarray*}

From the above, the result follows.
\end{Proof}

By the definition of the generalized Petersen graph $P(n,t)$, where $n \geq 3$ and $1 \leq t \leq \frac{n}{2}$, the following result is a straightforward observation:
\begin{itemize}
\item $P(n,t)\cong P(n,n-t)$;
\item $P(n,t)$ can be decomposed into three subgraphs with disjoint edges:
\begin{eqnarray*}
P(n,t)=G[V_1]\bigcup G[M]\bigcup G[V_{2}],
\end{eqnarray*}
\end{itemize}
where $G[V_1]$ is an $n$-cycle, $G[M]$ is the induced subgraph of a perfect matching $M$ in $P(n,t)$, $G[V_2]$ is composed by $(n,t)$ disjoint $\frac{n}{(n,t)}$-cycles, and $G[V_1]$ and $G[V_2]$ are edge-disjoint.

\begin{Theorem}
The generalized Petersen graph $G=P(n,t)$ is $(1,3)$-choosable.
\end{Theorem}
\begin{Proof}
Let $G$ be a graph with vertex set $$ V(G) = \{u_{1},u_{2},\ldots,u_{n}; v_{1},v_{2},\ldots,v_{n} \}$$
and edge set $$E(G) = \{u_{i}u_{i+1},u_{i}v_{i},v_{i}v_{i+t}: i \in \{1,2,\ldots,n \}\},$$ where the subscripts are taken modulo $n$.
Denote $V_{1}=\{u_{1},u_{2},\ldots,u_{n}\}$ as outer vertices set and $V_{2}=\{v_{1},v_{2},\ldots,v_{n}\}$ as inner vertices set. The edges $u_{i}u_{i+1},v_{i}v_{i+t}$, and $u_{i}v_{i}$ are denoted as outer edges, inner edges and leg edges, respectively, where $i \in \{1,2,\ldots, n\}$.

We will consider a demonstration of $P(n,t)$ into several types of cases in the following. Note that $n\geq 2t$.
If $n=2t$ and $t\ge 2$, then $\mad(G)\le \frac{5}{2}<\frac{11}{4}$.
By Theorem~\ref{mad}, $P(2t,t)$ is $(1,3)$-choosable.

We may assume $n\geq 2t+1$. Let $G'=G[E(G)\setminus E(H)]$.

If $t=1$, then choose  $H=G[V(1)\bigcup\{v_{1},v_{2}\}]$. By lemma~\ref{L1}$(i)$, $H$ is a $(1,3)$-choosable induced subgraph in $P(n,1)$. According to the structure of $P(n,1)$, it is not difficult to find $G'$ is a tree obtained by attaching an edge to each inner vertex of a path $P_{n-1}=v_2\cdots v_{n}v_1$. Then $G'$ is $(1,3)$-choosable by Lemma~\ref{tree}. Thus by Lemma~\ref{key1}, $P(n,1)$ is $(1,3)$-choosable. Take $P(5,1)$ as an example, see Figure $2$.
\begin{figure}[htbp]
\centering
\begin{minipage}[t]{0.24\textwidth}
\centering
\includegraphics[scale=0.6]{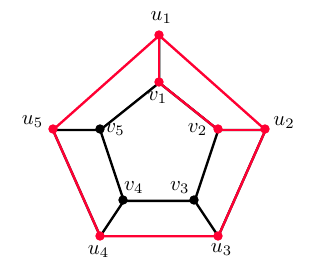}\\
\small{Fig. 2. $G=P(5,1)$}
\end{minipage}
\begin{minipage}[t]{0.24\textwidth}
\centering
\includegraphics[scale=0.9]{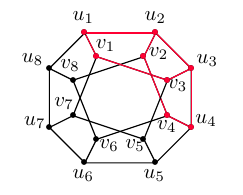}\\
\small{Fig. 3. $G=P(8,2)$}
\end{minipage}
\begin{minipage}[t]{0.23\textwidth}
\centering
\includegraphics[scale=0.6]{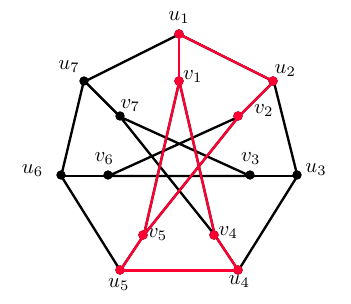}\\
\small{Fig. 4. $G=P(7,3)$}
\end{minipage}
\begin{minipage}[t]{0.25\textwidth}
\centering
\includegraphics[scale=0.6]{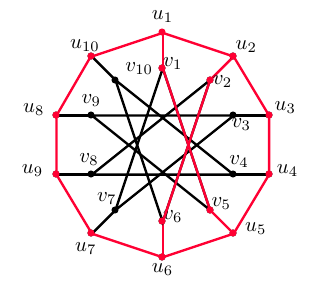}\\
\small{Fig. 5. $G=P(10,4)$}
\end{minipage}
\end{figure}

If $t=2$, then choose $H=G[X]$, where $X=\{u_{i},v_{i}:i=1,2,3,4\}$ and the edge set $E[X]$. Expanding the plane of the structure $H$, it is not difficult to see that $H$ is a graph obtained by two adjacent $5$-cycles with a common edge $u_2u_3$, by Lemma~\ref{L1}$(ii)$, $H$ is a $(1,3)$-choosable induced subgraph in $G$. According to the structure of $P(n,2)$, we claim that $G'$ is $(1,3)$-choosable. When $5=2t+1\leq n \leq 6$, $G'$ is a tree. When $n\geq7$, $G'$ contains an odd $5$-cycle $u_5u_6u_7v_7v_5u_5$, which is a non-bipartite $2$-degenerate graph. Thus by Theorems~\ref{tree} and~\ref{2d}, the claim is proved. By Lemma~\ref{key1}, $G=P(n,2)$ is $(1,3)$-choosable. Take $P(8,2)$ as an example, see Figure $3$.

Suppose that $t\geq3$.

If $n=2t+1$, then choose $H=G[X]$, where $X=\{u_{i},v_{i}:i=1,2,t+1,t+2\}$. Expanding the plane of the structure $H$, it is not difficult to see that $H$ is a graph obtained by two adjacent  $5$-cycles with a common edge $v_1v_{t+2}$. By Lemma~\ref{L1}$(ii)$, $H$ is a $(1,3)$-choosable induced subgraph in $G$. Note that $G'$ is a non-bipartite $2$-degenerate graph, since it contains an odd $5$-cycle $\{v_{3},v_{t+3},u_{t+3},u_{t+4},v_{t+4}v_3\}$. Then by Theorem~\ref{2d}, $G'$ is $(1,3)$-choosable. Therefore by Lemma~\ref{key1}, $G=P(2t+1,t)$  is $(1,3)$-choosable. Take $P(7,3)$ as an example, see Figure $4$.

If $n\geq 2t+2$, then choose $H=G[X]$, where $X=V_{1}\bigcup\{v_{i}:i=1,2,t+1,t+2\}$. Expanding the plane of the structure $H$, it is not difficult to see that it is isomorphic to the graph as in Lemma~\ref{L1}$(iii)$. Thus $H$ is a $(1,3)$-choosable induced subgraph in $G$. Note that $G[V_2]$ is a graph composed by $(n,t)$ disjoint $\frac{n}{(n,t)}$-cycles $C(1),\cdots,C((n,t))$. Denote by $C^*(i)$ be a unicyclic graph obtained by attaching one edge to each vertex of $\frac{n}{(n,t)}$-cycle $C(i)$. Then $G\setminus G[V_1]$ is a graph composed by $(n,t)$ disjoint $\frac{n}{(n,t)}$-cycles $C^*(1),\cdots,C^*((n,t))$, in which each $C^*(i)$ contains $v_i$ and $v_{i+t}$ for $i\in\{1,\dots,(n,t)\}$. According to the structure of $$G'=(G[M]-\{u_iv_i: i\in\{1,2,t+1,t+2\}\})\bigcup(G[V_2]-\{v_{1}v_{t+1}, v_{2}v_{t+2}\}),$$ we have three cases to discuss.
When $(n,t)=1$, $G'$ is a forest obtained from $C^*(1)$ by deleting two edges $v_{1}v_{t+1}$ and $v_{2}v_{t+2}$ in the $n$-cycle and four 
leg edges $\{u_iv_i: i\in\{1,2,t+1,t+2\}\}$ in $M$.
When $(n,t)=2$, $G'$ is a forest obtained from two disjoint unicyclic graphs $C^*(i)$ ($i=1,2$) by deleting one $v_{i}v_{t+i}$ in the $\frac{n}{2}$-cycle and two 
leg edges $\{u_iv_i: i\in\{i,t+i\}\}$ in $M$, respectively. Otherwise, $G'=F_0\cup G_0$, where $F_0$ is a forest consisting of two disjoint unicyclic graphs $C^*(i)$ ($i=1,2$) by deleting one edge $v_{i}v_{t+i}$ in the $\frac{n}{(n,t)}$-cycle and two 
leg edges $\{u_iv_i: i\in\{i,t+i\}\}$ in $M$, respectively,  and $G_0$ is a graph obtained from  $(n,t)-2$ disjoint unicyclic graphs $C^*(3),\cdots,C^*((n,t))$.
Note that $\mad(G')\le 2$. Then $G'$ is $(1,3)$-choosable by  Theorem~\ref{mad}.
Thus by Lemma~\ref{key1}, the result follows. Take $P(10,4)$ as an example, see Figure $5$.

From the above, the result follows.
\end{Proof}

\end{document}